# STATISTICAL INFERENCES FOR FUNCTIONAL DATA

By Jin-Ting Zhang[1] and Jianwei Chen

*National University of Singapore and University of Rochester*

With modern technology development, functional data are being observed frequently in many scientific fields. A popular method for analyzing such functional data is "smoothing first, then estimation." That is, statistical inference such as estimation and hypothesis testing about functional data is conducted based on the substitution of the underlying individual functions by their reconstructions obtained by one smoothing technique or another. However, little is known about this substitution effect on functional data analysis. In this paper this problem is investigated when the local polynomial kernel (LPK) smoothing technique is used for individual function reconstructions. We find that under some mild conditions, the substitution effect can be ignored asymptotically. Based on this, we construct LPK reconstruction-based estimators for the mean, covariance and noise variance functions of a functional data set and derive their asymptotics. We also propose a GCV rule for selecting good bandwidths for the LPK reconstructions. When the mean function also depends on some time-independent covariates, we consider a functional linear model where the mean function is linearly related to the covariates but the covariate effects are functions of time. The LPK reconstruction-based estimators for the covariate effects and the covariance function are also constructed and their asymptotics are derived. Moreover, we propose a $L^2$-norm-based global test statistic for a general hypothesis testing problem about the covariate effects and derive its asymptotic random expression. The effect of the bandwidths selected by the proposed GCV rule on the accuracy of the LPK reconstructions and the mean function estimator is investigated via a simulation study. The proposed methodologies are illustrated via an application to a real functional data set collected in climatology.

Received April 2004; revised May 2006.

[1]Supported by the National University of Singapore Academic Research Grant R-155-000-038-112.

**1. Introduction.** Functional data consist of functions which are often smooth but usually corrupted with noise. With modern technology development, such functional data are being observed frequently in many scientific fields; see Besse and Ramsay [3], Ramsay [20] and Ramsay and Dalzell [21], among others, for good examples and analyses. Comprehensive surveys about functional data analysis (FDA) can be found in [23, 24].

Mathematically, the above-mentioned functional data may be modeled as independent realizations of an underlying stochastic process,

$$(1.1) \qquad y_i(t) = \eta(t) + v_i(t) + \varepsilon_i(t), \qquad i = 1, 2, \ldots, n,$$

where $\eta(t)$ models the population mean function of the stochastic process, $v_i(t)$ is the $i$th individual variation (subject-effect) from $\eta(t)$, $\varepsilon_i(t)$ is the $i$th measurement error process and $y_i(t)$ is the $i$th response process. Without loss of generality, throughout this paper we assume the stochastic process has finite support, that is, $t \in \mathcal{T} = [a, b]$, $-\infty < a < b < \infty$. Moreover, we assume $v_i(t)$ and $\varepsilon_i(t)$ are independent, and are independent copies of $v(t) \sim \mathrm{SP}(0, \gamma)$ and $\varepsilon(t) \sim \mathrm{SP}(0, \gamma_\varepsilon), \gamma_\varepsilon(s, t) = \sigma^2(t) 1_{\{s=t\}}$, respectively, where and throughout $\mathrm{SP}(\eta, \gamma)$ denotes a stochastic process with mean function $\eta(t)$ and covariance function $\gamma(s, t)$. It follows that the underlying individual functions (trajectories) $f_i(t) = \mathrm{E}\{y_i(t)|v_i(t)\} = \eta(t) + v_i(t)$ are i.i.d. copies of the underlying stochastic process, $f(t) = \eta(t) + v(t) \sim \mathrm{SP}(\eta, \gamma)$. In practice, functional data are observed discretely. Let $t_{ij}, j = 1, 2, \ldots, n_i$, be the design time points of the $i$th subject. Then by (1.1) and letting $y_{ij} = y_i(t_{ij})$ and $\varepsilon_{ij} = \varepsilon_i(t_{ij})$, we have

$$(1.2) \quad y_{ij} = \eta(t_{ij}) + v_i(t_{ij}) + \varepsilon_{ij}, \qquad j = 1, 2, \ldots, n_i;\ i = 1, 2, \ldots, n.$$

In many practical situations, the above discrete functional data (1.2) have to be first registered before any statistical inference can be conducted. Methods for curve registration can be found in Kneip and Gasser [19], Kneip and Engel [18], Silverman [26], Ramsay and Silverman ([23], Chapter 5), Ramsay and Li [22] and Ramsay and Silverman ([24], Chapter 7), among others. In this paper, for convenience, we assume that the functional data (1.2) do not need registration or have been registered.

Estimation of the population characteristics $\eta(t), \gamma(s, t)$ and $\sigma^2(t)$ of the model (1.1) has been the focus of FDA in the literature. Most of the existing approaches involve one smoothing method or another. For example, Besse and Ramsay [3], Ramsay [20] and Ramsay and Dalzell [21] made use of reproducing kernel Hilbert space decomposition; Rice and Silverman [25] and Brumback and Rice [4] employed smoothing splines; Besse, Cardot and Ferraty [2] used B-splines; Hart and Wehrly [16] employed kernel smoothing; and Kneip [17] studied a principal components-based approach. Development of significance tests about $\eta(t)$ and other population characteristics of the model (1.1) is more important and challenging. Faraway [13] discussed the difficulties in extending some multivariate hypothesis testing procedures to FDA. Ramsay and Silverman [23] suggested a pointwise $t$-test or $F$-test but they did not discuss global tests. For curve data from stationary Gaussian processes, Fan and Lin [11] developed an adaptive Neyman test.

In this paper we adopt the method of "smoothing first, then estimation" for functional data. That is, we construct the estimators for $\eta(t), \gamma(s,t)$ and $\sigma^2(t)$ using the reconstructed individual functions $\hat{f}_i(t), i = 1, 2, \ldots, n$, obtained using one smoothing method or another; in particular, in this paper we use the local polynomial kernel (LPK) smoothing technique as described in [10], among others. The idea of "smoothing first, then estimation" itself is hardly new since it has been used in the literature; see [23, 24] and the references therein. What is new here is that we investigate the effect of the substitution of the underlying individual functions $f_i(t), i = 1, 2, \ldots, n$, by their LPK reconstructions in FDA. We show that, under some mild conditions, the effect of such a substitution is asymptotically ignorable in FDA. Based on this, we derive the asymptotics of the estimators $\hat{\eta}(t), \hat{\gamma}(s,t)$ and $\hat{\sigma}^2(t)$. In particular, under some mild conditions, we show that: (1) $\hat{\eta}(t)$ and $\hat{\gamma}(s,t)$ are $\sqrt{n}$-consistent and asymptotically Gaussian; (2) the asymptotic efficiency of $\hat{\eta}(t)$ will not be affected by better choice of the bandwidth than the bandwidth selected by a GCV rule; and (3) the convergence rate of $\hat{\sigma}^2(t)$ is affected by the convergence rate of the LPK reconstructions. More details about these results are given in Section 2.

In the model (1.1) the only covariate for the mean function $\eta(t)$ is time. In many applications $\eta(t)$ may also depend on some time-independent covariates and can be written as $\eta(t; \mathbf{x}) = \mathbf{x}^T \boldsymbol{\beta}(t)$, where the covariate vector $\mathbf{x} = [x_1, \ldots, x_q]^T$ and the unknown but smooth coefficient function vector $\boldsymbol{\beta}(t) = [\beta_1(t), \ldots, \beta_q(t)]^T$. A replacement of $\eta(t)$ by $\eta(t; \mathbf{x}_i) = \mathbf{x}_i^T \boldsymbol{\beta}(t)$ in (1.1) leads to the so-called functional linear model

$$(1.3) \qquad y_i(t) = \mathbf{x}_i^T \boldsymbol{\beta}(t) + v_i(t) + \varepsilon_i(t), \qquad i = 1, 2, \ldots, n,$$

where $y_i(t), v_i(t)$ and $\varepsilon_i(t)$ are the same as those defined in (1.1). The ignorability of the substitution effect is also applied to the LPK reconstructions $\hat{f}_i(t)$ of the individual functions $f_i(t) = \mathbf{x}_i^T \boldsymbol{\beta}(t) + v_i(t)$ of the above model. Based on this, we construct the estimators $\hat{\boldsymbol{\beta}}(t)$ and $\hat{\gamma}(s,t)$ and investigate their asymptotics; in particular, we show that $\hat{\boldsymbol{\beta}}(t)$ is $\sqrt{n}$-consistent and asymptotically Gaussian. Moreover, we propose a global $L^2$-norm-based test statistic $T_n$ to test a general hypothesis testing problem about the covariate effects $\boldsymbol{\beta}(t)$; its asymptotic random expression is derived. More details about these results are given in Section 3.

The rest of the paper is organized as follows. In Section 4 we present a simulation study which aims to investigate the effect of the bandwidth choice on the accuracy of the LPK reconstructions $\hat{f}_i(t)$ and the mean function estimator $\hat{\eta}(t)$. In Section 5 we illustrate the proposed methodologies by applying them to a real functional data set collected in climatology. Finally, in Section 6 technical proofs of some asymptotic results are outlined.

## 2. Basic methodologies.

2.1. *LPK reconstruction of individual functions.* First of all, we describe how to reconstruct the individual functions $f_i(t), i = 1, 2, \ldots, n$, using the LPK smoothing technique based on the standard nonparametric regression model

$$(2.1) \qquad y_{ij} = f_i(t_{ij}) + \varepsilon_{ij}, \qquad j = 1, 2, \ldots, n_i;\ i = 1, 2, \ldots, n.$$

For any fixed time point $t$, assume $f_i(t)$ has a $(p+1)$th continuous derivative in a neighborhood of $t$ for some positive integer $p$. Then by Taylor's expansion, $f_i(t_{ij})$ can be locally approximated by a $p$-order polynomial, that is,

$$f_i(t_{ij}) \approx f_i(t) + (t_{ij} - t) f_i^{(1)}(t) + \cdots + (t_{ij} - t)^p f_i^{(p)}(t)/p! = \mathbf{z}_{ij}^T \boldsymbol{\alpha}_i,$$

in the neighborhood of $t$, where $\boldsymbol{\alpha}_i = [\alpha_{i0}, \alpha_{i1}, \ldots, \alpha_{ip}]^T$ with $\alpha_{ir} = f_i^{(r)}(t)/r!$, and $\mathbf{z}_{ij} = [1, t_{ij} - t, \ldots, (t_{ij} - t)^p]^T$. Then the $p$-order LPK reconstructions of $f_i(t)$ are defined as $\hat{f}_i(t) = \hat{\alpha}_{i0} = \mathbf{e}_{1,p+1}^T \hat{\boldsymbol{\alpha}}_i$, where and throughout $\mathbf{e}_{r,s}$ denotes the $s$-dimensional unit vector whose $r$th component is 1 and others are 0, and $\hat{\boldsymbol{\alpha}}_i$ are the minimizers of the weighted least squares criterion

$$(2.2) \qquad \begin{aligned} &\sum_{i=1}^{n} \sum_{j=1}^{n_i} [y_{ij} - \mathbf{z}_{ij}^T \boldsymbol{\alpha}_i]^2 K_h(t_{ij} - t) \\ &= \sum_{i=1}^{n} (\mathbf{y}_i - \mathbf{Z}_i \boldsymbol{\alpha}_i)^T \mathbf{K}_{ih} (\mathbf{y}_i - \mathbf{Z}_i \boldsymbol{\alpha}_i), \end{aligned}$$

where $\mathbf{y}_i = [y_{i1}, \ldots, y_{in_i}]^T, \mathbf{Z}_i = [\mathbf{z}_{i1}, \ldots, \mathbf{z}_{in_i}]^T$ and $\mathbf{K}_{ih} = \operatorname{diag}(K_h(t_{i1} - t), \ldots, K_h(t_{in_i} - t))$, with $K_h(\cdot) = K(\cdot/h)/h$, obtained by rescaling a kernel function $K(\cdot)$ (often a symmetric p.d.f.) with bandwidth $h > 0$ that controls the size of the associated neighborhood. Minimizing (2.2) with respect to $\boldsymbol{\alpha}_i, i = 1, 2, \ldots, n$, is equivalent to minimizing the $i$th term in the summation on the right-hand side of (2.2) with respect to $\boldsymbol{\alpha}_i$ for each $i = 1, 2, \ldots, n$. It follows that for $i = 1, 2, \ldots, n$,

$$(2.3) \qquad \begin{aligned} \hat{f}_i(t) &= \mathbf{e}_{1,p+1}^T (\mathbf{Z}_i^T \mathbf{K}_{ih} \mathbf{Z}_i)^{-1} \mathbf{Z}_i^T \mathbf{K}_{ih} \mathbf{y}_i \\ &= \sum_{j=1}^{n_i} K_h^{n_i}(t_{ij} - t) y_{ij}, \end{aligned}$$

where $K^{n_i}(t)$ are known as the empirical equivalent kernels for the $p$-order LPK; see Fan and Gijbels [10].

In (2.2) different bandwidths may be used for different individual functions. However, the individual functions in a functional data set are i.i.d. realizations of a stochastic process, and hence, often admit similar smoothness properties and sometimes similar shapes [17, 18]; it is then reasonable to treat them in the same way, for example, using a common bandwidth for all of them. The advantages in using a common bandwidth at least include the following: (a) reduce the computational effort for bandwidth selection; and (b) simplify the asymptotic results of the estimators.

For convenience, we define the following widely-used functionals of a kernel $K$:

$$B_r(K) = \int K(t) t^r \, dt,$$
$$(2.4) \qquad V(K) = \int K(t)^2 \, dt,$$
$$K^{(1)}(t) = \int K(s) K(s+t) \, ds.$$

For estimating a function instead of derivatives, Fan and Gijbels [10] pointed out that even orders are not appealing. Therefore, throughout this paper, we assume $p$ is an odd integer; moreover, we denote $\gamma_{k,l}(s,t)$ as the $(k,l)$-times partial derivative of $\gamma(s,t)$, that is, $\gamma_{k,l}(s,t) = \frac{\partial^{k+l}\gamma(s,t)}{\partial^k s \, \partial^l t}$, and denote $\mathcal{D}$ as the set of all the design time points $t_{ij}, j = 1, 2, \ldots, n_i;\ i = 1, 2, \ldots, n$. In addition, we denote $O_{\mathrm{UP}}(1)$ [resp. $o_{\mathrm{UP}}(1)$] as "bounded (resp. tends to 0) in probability uniformly for any $t$ within the interior of $\mathcal{T}$ and all $i = 1, 2, \ldots, n$." Finally, the following regular conditions are imposed.

CONDITION A.

1. The design time points $t_{ij}, j = 1, 2, \ldots, n_i;\ i = 1, 2, \ldots, n$, are i.i.d. with p.d.f. $\pi(\cdot)$ which has the bounded support $\mathcal{T} = [a, b]$. For any given $t$ within the interior of $\mathcal{T}$, $\pi'(t)$ exists and is continuous over $\mathcal{T}$.
2. Let $s$ and $t$ be any two interior time points of $\mathcal{T}$. The individual functions $f_i(t), i = 1, 2, \ldots, n$, and their mean function $\eta(t)$ have up to $(p+1)$-times continuous derivatives. Their covariance function $\gamma(s,t)$ has up to $(p+1)$-times continuous derivatives for both $s$ and $t$. The variance function of the measurement errors, $\sigma^2(t)$, is continuous at $t$.
3. The kernel $K$ is a bounded symmetrical p.d.f. with bounded support $[-1, 1]$.
4. There are two positive constants $C$ and $\delta$ such that $n_i \geq C n^\delta$, for all $i = 1, 2, \ldots, n$. As $n \to \infty$, we have $h \to 0$ and $n^\delta h \to \infty$.

REMARK 1. For some practical functional data sets, Condition A4 may be too restrictive. For example, a functional data set with a few individual functions having $n_i < Cn^\delta$ does not satisfy Condition A4. However, such a functional data set can often be slightly modified to satisfy Condition A4. A simple way of doing so is to drop those individual functions having $n_i < Cn^\delta$ so that the remaining individual functions form a new functional data set which satisfies Condition A4. This procedure will not result in less efficient estimators when $\tilde{n}/n \to 0$ and will not affect the consistency of the estimators when $(n - \tilde{n}) \to \infty$, where $\tilde{n}$ is the number of dropped individual functions, which may be bounded or tend to $\infty$ as $n \to \infty$.

Using Lemma 3 in Section 6, it is easy to show the following.

THEOREM 1. *Assume Condition A is satisfied. Then the average conditional MSE (mean squared errors) of the p-order LPK reconstructions $\hat{f}_i(t), i = 1, 2, \ldots, n$, is*

$$
\begin{aligned}
& \mathrm{E}\left\{n^{-1} \sum_{i=1}^{n} [\hat{f}_i(t) - f_i(t)]^2 \Big| \mathcal{D}\right\} \\
(2.5) \quad & = \left\{ \frac{B_{p+1}^2(K^*)[(\eta^{(p+1)}(t))^2 + \gamma_{p+1,p+1}(t,t)]}{(p+1)!^2} h^{2(p+1)} \right. \\
& \qquad \left. + \frac{V(K^*)\sigma^2(t)}{\pi(t)} (\tilde{m}h)^{-1} \right\} [1 + o_P(1)],
\end{aligned}
$$

*where $K^*$ is the equivalent kernel of the p-order LPK ([10], page 64), and $\tilde{m} = (n^{-1} \sum_{i=1}^{n} n_i^{-1})^{-1}$.*

REMARK 2. On the left-hand side of (2.5) the notation $\mathrm{E}\{\cdot|\mathcal{D}\}$ denotes the conditional expectation when all the design time points $t_{ij}, j = 1, 2, \ldots, n_i;\ i = 1, 2, \ldots, n$, are given. Nevertheless, the leading term on the right-hand side of (2.5) is independent of $\mathcal{D}$ and hence, the left-hand side is nearly unconditional. For technical convenience and following the literature tradition (e.g., [10]), we keep using the "conditional expectation" notation $\mathrm{E}\{\cdot|\mathcal{D}\}$ here and throughout. This remark applies to all other statistical operations conditional to $\mathcal{D}$ given in this paper.

Theorem 1 indicates that the optimal bandwidth of the $p$-order LPK reconstructions $\hat{f}_i(t)$ is $h = O_P(\tilde{m}^{-1/(2p+3)}) = O_P(n^{-\delta/(2p+3)})$. Using Lemma 3 again, we can show the following.

THEOREM 2. *Assume Condition A is satisfied. Then for the $p$-order LPK reconstructions $\hat{f}_i(t)$ using the bandwidth $h = O(n^{-\delta/(2p+3)})$, we have*

$$(2.6) \qquad \hat{f}_i(t) = f_i(t) + n^{-(p+1)\delta/(2p+3)} O_{\mathrm{UP}}(1), \qquad i = 1, 2, \ldots, n.$$

Theorem 2 implies that, under the given conditions, the LPK reconstructions $\hat{f}_i(t)$ are asymptotically uniformly little different from the underlying individual functions $f_i(t)$. We expect that this is true not only for LPK but also for any other linear smoothers, for example, smoothing splines [14, 27], regression splines or orthogonal series [7], among others.

2.2. *Estimation of the mean and covariance functions.* It is then natural to estimate the mean function $\eta(t)$ and the covariance function $\gamma(s,t)$ by the sample mean and sample covariance functions of the $p$-order LPK reconstructions $\hat{f}_i(t)$,

$$\hat{\eta}(t) = n^{-1} \sum_{i=1}^{n} \hat{f}_i(t),$$

(2.7)

$$\hat{\gamma}(s,t) = (n-1)^{-1} \sum_{i=1}^{n} \{\hat{f}_i(s) - \hat{\eta}(s)\}\{\hat{f}_i(t) - \hat{\eta}(t)\}.$$

The asymptotic conditional bias, covariance and variance for $\hat{\eta}(t)$ are given below.

THEOREM 3. *Assume Condition A is satisfied. Then as $n \to \infty$, the asymptotic conditional bias, covariance and variance of $\hat{\eta}(t)$ are*

$\text{Bias}\{\hat{\eta}(t)|\mathcal{D}\}$

$$= \frac{B_{p+1}(K^*)\eta^{(p+1)}(t)}{(p+1)!} h^{p+1}[1 + o_P(1)],$$

$\text{Cov}\{\hat{\eta}(s), \hat{\eta}(t)|\mathcal{D}\}$

$$= \left\{\gamma(s,t)/n + \frac{K^{*(1)}[(s-t)/h]\sigma^2(s)}{\pi(t)}(n\tilde{m}h)^{-1}\right.$$

$$\left. + \frac{B_{p+1}(K^*)[\gamma_{p+1,0}(s,t) + \gamma_{0,p+1}(s,t)]}{(p+1)!} n^{-1} h^{p+1}\right\}[1 + o_P(1)],$$

$\text{Var}\{\hat{\eta}(t)|\mathcal{D}\}$

$$= \left\{\gamma(t,t)/n + \frac{V(K^*)\sigma^2(t)}{\pi(t)}(n\tilde{m}h)^{-1}\right.$$

$$\left. + \frac{B_{p+1}(K^*)[\gamma_{p+1,0}(t,t) + \gamma_{0,p+1}(t,t)]}{(p+1)!} n^{-1} h^{p+1}\right\}[1 + o_P(1)].$$

REMARK 3. Under Condition A and by Theorem 3, we have

(2.8) $\text{MSE}\{\hat{\eta}(t)|\mathcal{D}\} = \gamma(t,t)/n + O_{\text{UP}}\{h^{2(p+1)} + (n\tilde{m}h)^{-1} + n^{-1}h^{p+1}\}.$

We then always have $\mathrm{MSE}\{\hat{\eta}(t)|\mathcal{D}\} = \gamma(t,t)/n + o_{\mathrm{UP}}(1/n)$, provided that

(2.9) $$\tilde{m}h \to \infty, \qquad nh^{2(p+1)} \to 0.$$

REMARK 4. Condition (2.9) is satisfied by any bandwidth $h = O(n^{-\delta^*})$ with $1/[2(p+1)] < \delta^* < \delta$. In particular, it is satisfied by the optimal bandwidth, $h = O(n^{-\delta/(2p+3)})$, for the $p$-order LPK reconstructions $\hat{f}_i(t)$ when $\delta > 1 + 1/[2(p+1)]$. In this case, the $p$-order LPK reconstruction optimal bandwidth is sufficiently small to guarantee the $\sqrt{n}$-consistency of $\hat{\eta}(t)$. Condition (2.9) is also satisfied by the optimal bandwidth, $h = O(n^{-(1+\delta)/(2p+3)})$ (when $1 + 1/[2(p+1)] < \delta < 1 + 1/(p+1)$) or $h = O(n^{-\delta/(p+2)})$ [when $\delta > 1 + 1/(p+1)$], for $\hat{\eta}(t)$. It follows that, in both cases, the optimal bandwidth admits the same asymptotic efficiency for estimating $\eta(t)$ because $\mathrm{MSE}\{\sqrt{n}\hat{\eta}(t)|\mathcal{D}\} \to \gamma(t,t)$ as $n \to \infty$.

By pretending all the underlying individual functions $f_i(t)$ were observed, the "ideal" estimators of $\eta(t)$ and $\gamma(s,t)$ are

(2.10)
$$\tilde{\eta}(t) = n^{-1} \sum_{i=1}^{n} f_i(t),$$

$$\tilde{\gamma}(s,t) = (n-1)^{-1} \sum_{i=1}^{n} \{f_i(s) - \tilde{\eta}(s)\}\{f_i(t) - \tilde{\eta}(t)\}.$$

THEOREM 4. *Assume Condition A is satisfied, and the bandwidth $h = O(n^{-\delta/(2p+3)})$ is used for the p-order LPK reconstructions $\hat{f}_i(t)$. Then as $n \to \infty$, we have*

(2.11)
$$\hat{\eta}(t) = \tilde{\eta}(t) + n^{-(p+1)\delta/(2p+3)} O_{\mathrm{UP}}(1),$$

$$\hat{\gamma}(s,t) = \tilde{\gamma}(s,t) + n^{-(p+1)\delta/(2p+3)} O_{\mathrm{UP}}(1).$$

*In addition, assume $\delta > 1 + 1/[2(p+1)]$. Then as $n \to \infty$, we have*

(2.12)
$$\sqrt{n}\{\hat{\eta}(t) - \eta(t)\} \sim \mathrm{AGP}(0, \gamma),$$

$$\sqrt{n}\{\hat{\gamma}(s,t) - \gamma(s,t)\} \sim \mathrm{AGP}(0, \gamma^*),$$

*where $\mathrm{AGP}(\eta, \gamma)$ denotes an asymptotic Gaussian process with mean function $\eta(t)$ and covariance function $\gamma(s,t)$, and*

(2.13) $\gamma^*\{(s_1, t_1), (s_2, t_2)\} = \mathrm{E}\{v_1(s_1)v_1(t_1)v_1(s_2)v_1(t_2)\} - \gamma(s_1, t_1)\gamma(s_2, t_2),$

*with $v_1(t)$ denoting the subject effect of the first individual function $f_1(t)$ as defined in (1.1). When the subject effect process $v(t)$ is Gaussian,*

$$\gamma^*\{(s_1, t_1), (s_2, t_2)\} = \gamma(s_1, t_2)\gamma(s_2, t_1) + \gamma(s_1, s_2)\gamma(t_1, t_2).$$

Theorem 4 indicates that, under some mild conditions, the proposed estimators (2.7) are asymptotically identical to the "ideal" estimators (2.10). The required key condition is $\delta > 1 + 1/[2(p+1)]$. It follows that, to make the measurement errors ignorable via LPK smoothing, we need the number of measurements, $n_i$, for all the subjects (or a large number of subjects; see Remark 1 for discussion) to tend to infinity slightly faster than the number of subjects, $n$.

2.3. *Estimation of the noise variance function.* The noise variance function $\sigma^2(t)$ measures the variation of the measurement errors $\varepsilon_{ij}$ of the model (1.2). Following Hall and Marron [15] and Fan and Yao [12], we can construct a $\tilde{p}$-order LPK estimator of $\sigma^2(t)$ based on the $p$-order LPK residuals $\hat{\varepsilon}_{ij} = y_{ij} - \hat{f}_i(t_{ij})$, although our setting is more complicated. As expected, the resulting $\tilde{p}$-order LPK estimator of $\sigma^2(t)$ will be consistent, but its convergence rate will be affected by that of the $p$-order LPK reconstructions $\hat{f}_i(t), i = 1, 2, \ldots, n$.

As an illustration, let us consider the simplest LPK estimator, that is, the kernel estimator for $\sigma^2(t)$ based on $\hat{\varepsilon}_{ij}$,

$$(2.14) \qquad \hat{\sigma}^2(t) = \frac{\sum_{i=1}^n \sum_{j=1}^{n_i} H_b(t_{ij} - t) \hat{\varepsilon}_{ij}^2}{\sum_{i=1}^n \sum_{j=1}^{n_i} H_b(t_{ij} - t)},$$

where $H_b(\cdot) = H(\cdot/b)/b$ with the kernel function $H$ and the bandwidth $b$.

Pretending $\hat{\varepsilon}_{ij} \equiv \varepsilon_{ij}$, by standard kernel estimation theory (Wand and Jones [28] and Fan and Gijbels [10], among others), the optimal bandwidth for $\hat{\sigma}^2(t)$ is $b = O_P(N^{-1/5})$, where $N = \sum_{i=1}^n n_i$ denotes the total number of measurements for all the subjects, and the associated convergence rate of $\hat{\sigma}^2(t)$ is $O_P(N^{-2/5})$. However, for the current setup, this convergence rate will be affected by the convergence rate of the $p$-order LPK reconstructions $\hat{f}_i(t), i = 1, 2, \ldots, n$, since under Condition A and by Theorem 2, we actually only have $\hat{\varepsilon}_{ij} = \varepsilon_{ij} + n^{-(p+1)\delta/(2p+3)} O_{\mathrm{UP}}(1)$. For convenience, let $\nu_1(t) = \mathrm{E}[\varepsilon_i^2(t)] = \sigma^2(t)$ and $\nu_2(t) = \mathrm{Var}[\varepsilon_i^2(t)]$.

THEOREM 5. *Assume Condition A is satisfied and the p-order LPK reconstructions $\hat{f}_i(t)$ use a bandwidth $h = O(n^{-\delta/(2p+3)})$. In addition, assume $\nu_1'(t)$ and $\nu_2(t)$ exist and are continuous at $t \in \mathcal{T}$, and the kernel estimator $\hat{\sigma}^2(t)$ uses a bandwidth $b = O(N^{-1/5})$. Then we have*

$$(2.15) \qquad \hat{\sigma}^2(t) = \sigma^2(t) + O_{\mathrm{UP}}(n^{-2(1+\delta)/5} + n^{-(p+1)\delta/(2p+3)}).$$

By the above theorem, it is seen that when $\delta < 2(2p+3)/(p-1)$, the second order term dominates the first order term; and in particular, when $p = 1$, we have $\hat{\sigma}^2(t) = \sigma^2(t) + O_{\mathrm{UP}}(n^{-2\delta/5})$. In this case the optimal convergence rate of $\hat{\sigma}^2(t)$ is not attainable. It is attainable only when $\delta > 2(2p+3)/(p-1)$, so that the first order term in (2.15) dominates the second order term. This is the case only when $p \geq 3$. When $p = 3$, $\delta > 9$ is required; and when $p = 2k + 1 \to \infty$, $\delta > 4$ is required. Therefore, it is usually difficult to make the convergence rate of $\hat{\sigma}^2(t)$ unaffected by the convergence rate of the $p$-order LPK reconstructions $\hat{f}_i(t)$.

2.4. *Bandwidth selection.* Theorem 1 suggests that we can choose a good bandwidth for the $p$-order LPK reconstructions $\hat{f}_i(t), i = 1, 2, \ldots, n$, using the generalized cross-validation (GCV) score

$$\text{(2.16)} \qquad \text{GCV}(h) = n^{-1} \sum_{i=1}^{n} \text{GCV}_i(h),$$

where $\text{GCV}_i(h)$ is the GCV score of the $i$th $p$-order LPK reconstruction $\hat{f}_i(t)$. Let $\mathbf{A}_i$ be the smoother matrix of the $i$th subject constructed using (2.3). Then we have $\hat{\mathbf{y}}_i = \mathbf{A}_i \mathbf{y}_i$ and $\text{GCV}_i(h) = \mathbf{y}_i^T (\mathbf{I}_{n_i} - \mathbf{A}_i)^T (\mathbf{I}_{n_i} - \mathbf{A}_i) \mathbf{y}_i / [1 - \text{tr}(\mathbf{A}_i)/n_i]^2$, where $\mathbf{y}_i = [y_{i1}, \ldots, y_{in_i}]^T$, $\hat{\mathbf{y}}_i = [\hat{y}_{i1}, \ldots, \hat{y}_{in_i}]^T$ and $\text{tr}(\mathbf{S})$ denotes the trace of the matrix $\mathbf{S}$. In practice, the optimal bandwidth $h^*$ can be obtained via minimizing $\text{GCV}(h)$ over a number of bandwidth candidates of interest. Theoretically, it is expected that $h^* = O_P(n^{-\delta/(2p+3)})$.

Remark 4 states that, under the required conditions, the optimal bandwidth for $\hat{\eta}(t)$ and the optimal bandwidth for the $p$-order LPK reconstructions $\hat{f}_i(t)$ admit the same asymptotic efficiency for estimating $\eta(t)$. Therefore, it is generally sufficient to use $h^*$ for estimating $\eta(t)$ although, for finite samples, better bandwidth choices for $\hat{\eta}(t)$ are possible.

**3. Functional linear models.** Notice that Theorem 2 is also applied to the $p$-order LPK reconstructions $\hat{f}_i(t)$ of the underlying individual functions $f_i(t) = \mathbf{x}_i^T \boldsymbol{\beta}(t) + v_i(t), i = 1, 2, \ldots, n$, of the functional linear model (1.3). This property can be used to do inference about the model (1.3). In this section we focus on the estimation and significance tests of the coefficient function vector (covariate effects) $\boldsymbol{\beta}(t)$ of the model.

3.1. *Coefficient function estimation.* Let $\hat{\mathbf{f}}(t) = [\hat{f}_1(t), \ldots, \hat{f}_n(t)]^T$ and $\mathbf{X} = [\mathbf{x}_1, \ldots, \mathbf{x}_n]^T$. Throughout this paper we assume $\mathbf{X}$ has full rank. Then the least-squares estimator of $\boldsymbol{\beta}(t)$ is

$$\text{(3.1)} \qquad \hat{\boldsymbol{\beta}}(t) = \left\{ \sum_{i=1}^{n} \mathbf{x}_i \mathbf{x}_i^T \right\}^{-1} \sum_{i=1}^{n} \mathbf{x}_i \hat{f}_i(t) = (\mathbf{X}^T \mathbf{X})^{-1} \mathbf{X}^T \hat{\mathbf{f}}(t),$$

which minimizes $Q(\boldsymbol{\beta}) = n^{-1} \sum_{i=1}^{n} \int [\hat{f}_i(t) - \mathbf{x}_i^T \boldsymbol{\beta}(t)]^2 \, dt$. It follows that the subject-effects $v_i(t)$ can be estimated by $\hat{v}_i(t) = \hat{f}_i(t) - \mathbf{x}_i^T \hat{\boldsymbol{\beta}}(t)$ and their covariance function $\gamma(s,t)$ can be estimated by

$$\hat{\gamma}(s,t) = (n-q)^{-1} \sum_{i=1}^{n} \hat{v}_i(s)\hat{v}_i(t)$$

(3.2)
$$= (n-q)^{-1}\hat{\mathbf{v}}(s)^T\hat{\mathbf{v}}(t),$$

where $\hat{\mathbf{v}}(t) = [\hat{v}_1(t), \hat{v}_2(t), \ldots, \hat{v}_n(t)]^T = \hat{\mathbf{f}}(t) - \mathbf{X}(\mathbf{X}^T\mathbf{X})^{-1}\mathbf{X}^T\hat{\mathbf{f}}(t) = (\mathbf{I}_n - \mathbf{P})\hat{\mathbf{f}}(t)$ and $\mathbf{P} = \mathbf{X}(\mathbf{X}^T\mathbf{X})^{-1}\mathbf{X}^T$ is a projection matrix with $\mathbf{P}^T = \mathbf{P}, \mathbf{P}^2 = \mathbf{P}$ and $\text{tr}(\mathbf{P}) = q$.

Pretending $f_i(t), i = 1, 2, \ldots, n$, is known, the "ideal" estimators of $\boldsymbol{\beta}(t)$ and $\gamma(s,t)$ are

(3.3)
$$\tilde{\boldsymbol{\beta}}(t) = (\mathbf{X}^T\mathbf{X})^{-1}\mathbf{X}^T\mathbf{f}(t),$$
$$\tilde{\gamma}(s,t) = (n-q)^{-1}\tilde{\mathbf{v}}(s)^T\tilde{\mathbf{v}}(t),$$

where $\mathbf{f}(t) = [f_1(t), \ldots, f_n(t)]^T$ and $\tilde{\mathbf{v}}(t) = (\mathbf{I}_n - \mathbf{P})\mathbf{f}(t)$. It is easy to show that $\mathrm{E}\tilde{\boldsymbol{\beta}}(t) = \boldsymbol{\beta}(t)$ and $\mathrm{E}\tilde{\gamma}(s,t) = \gamma(s,t)$. For further investigation, we impose the following conditions.

CONDITION B.

1. The covariate vectors $\mathbf{x}_i, i = 1, 2, \ldots, n$, are i.i.d. with finite and invertible second moment $\mathrm{E}\mathbf{x}_1\mathbf{x}_1^T = \boldsymbol{\Omega}$; moreover, they are uniformly bounded in probability; that is, $\mathbf{x}_i = O_{\mathrm{UP}}(1)$.
2. The subject-effects $v_i(t)$ are uniformly bounded in probability; that is, $\mathbf{v}_i(t) = O_{\mathrm{UP}}(1)$.

THEOREM 6. *Assume Conditions* A *and* B *are satisfied, and the p-order LPK reconstructions* $\hat{f}_i(t)$ *use a bandwidth* $h = O(n^{-\delta/(2p+3)})$. *Then as* $n \to \infty$, *we have*

(3.4)
$$\hat{\boldsymbol{\beta}}(t) = \tilde{\boldsymbol{\beta}}(t) + n^{-(p+1)\delta/(2p+3)}O_{\mathrm{UP}}(1),$$
$$\hat{\gamma}(s,t) = \tilde{\gamma}(s,t) + n^{-(p+1)\delta/(2p+3)}O_{\mathrm{UP}}(1).$$

*In addition, assume* $\delta > 1 + 1/[2(p+1)]$. *Then as* $n \to \infty$, *we have*

(3.5) $$\sqrt{n}\{\hat{\boldsymbol{\beta}}(t) - \boldsymbol{\beta}(t)\} \sim \mathrm{AGP}(0, \gamma_\beta),$$

*where* $\gamma_\beta(s,t) = \gamma(s,t)\boldsymbol{\Omega}^{-1}$.

Theorem 6 implies that, under the given conditions, the proposed estimators $\hat{\boldsymbol{\beta}}(t)$ and $\hat{\gamma}(s,t)$ are asymptotically identical to the "ideal" estimators $\tilde{\boldsymbol{\beta}}(t)$ and $\tilde{\gamma}(s,t)$, respectively. Therefore, in FDA it seems reasonable to directly assume the underlying individual functions are "observed" as is done in [23, 24]. The asymptotic result stated in (3.5) is a foundation for significance tests of the covariate effects.

3.2. *Significance tests of the covariate effects.* Consider the general hypothesis testing problem

(3.6) $$H_0 : \mathbf{C}\boldsymbol{\beta}(t) = \mathbf{c}(t), \quad \text{vs.} \quad H_1 : \mathbf{C}\boldsymbol{\beta}(t) \neq \mathbf{c}(t),$$

where $t \in \mathcal{T} = [a,b]$, $\mathbf{C}$ is a given $k \times q$ full rank matrix, and $\mathbf{c}(t) = [c_1(t), \ldots, c_k(t)]^T$ is a given vector of functions. In order to check the significance of the $r$th covariate effect, one takes $\mathbf{C} = \mathbf{e}_{r,q}^T = [0, \ldots, 0, 1, 0, \ldots, 0]$ and $\mathbf{c}(t) = 0$; in order to check if the first two coefficient functions are the same, that is, $\beta_1(t) = \beta_2(t)$, one takes $\mathbf{C} = (\mathbf{e}_{1,q} - \mathbf{e}_{2,q})^T = [1, -1, 0, \ldots, 0]$ and $\mathbf{c}(t) = 0$.

It is natural to estimate $\mathbf{C}\boldsymbol{\beta}(t)$ by $\mathbf{C}\hat{\boldsymbol{\beta}}(t)$. By Theorem 6, we have

(3.7) $$\sqrt{n}[\mathbf{C}\hat{\boldsymbol{\beta}}(t) - \mathbf{c}(t)] \sim \text{AGP}(\boldsymbol{\eta}_c, \boldsymbol{\gamma}_c),$$

where $\boldsymbol{\eta}_c(t) = \sqrt{n}[\mathbf{C}\boldsymbol{\beta}(t) - \mathbf{c}(t)]$ and $\boldsymbol{\gamma}_c(s,t) = \gamma(s,t)\mathbf{C}\boldsymbol{\Omega}^{-1}\mathbf{C}^T$. Let

(3.8) $$\begin{aligned}\mathbf{w}(t) &= \{\mathbf{C}(\mathbf{X}^T\mathbf{X})^{-1}\mathbf{C}^T\}^{-1/2}[\mathbf{C}\hat{\boldsymbol{\beta}}(t) - \mathbf{c}(t)] \\ &= [w_1(t), \ldots, w_k(t)]^T.\end{aligned}$$

Since $\mathbf{X}^T\mathbf{X}/n \to \boldsymbol{\Omega}$ as $n \to \infty$, using (3.7), we can show that $\mathbf{w}(t) \sim \text{AGP}(\boldsymbol{\eta}_w, \boldsymbol{\gamma}_w)$, where

(3.9) $$\begin{aligned}\boldsymbol{\eta}_w(t) &= \sqrt{n}(\mathbf{C}\boldsymbol{\Omega}^{-1}\mathbf{C}^T)^{-1/2}[\mathbf{C}\boldsymbol{\beta}(t) - \mathbf{c}(t)] \\ &= [\eta_{w1}(t), \ldots, \eta_{wk}(t)]^T, \\ \boldsymbol{\gamma}_w(s,t) &= \gamma(s,t)\mathbf{I}_k,\end{aligned}$$

where $\mathbf{I}_k$ denotes the identity matrix of size $k$. It follows that the components $w_1(t), \ldots, w_k(t)$ are independent asymptotic Gaussian processes with mean functions $\eta_{w1}(t), \ldots, \eta_{wk}(t)$, respectively, and a common covariance function $\gamma(s,t)$. That is,

(3.10) $$w_l(t) \sim \text{AGP}(\eta_{wl}, \gamma), \qquad l = 1, 2, \ldots, k.$$

Based on these results and with $\mathbf{C}$ and $\mathbf{c}(t)$ properly specified, pointwise $t$ and $F$-tests for the coefficient functions $\beta_1(t), \ldots, \beta_q(t)$ can easily be conducted ([23], Chapter 9). We here propose the following global test statistic for the general hypothesis testing problem (3.6):

(3.11) $$T_n = \int_a^b \|\mathbf{w}(t)\|^2 \, dt = \sum_{l=1}^k \int_a^b w_l^2(t) \, dt,$$

where $\|\cdot\|$ denotes the usual $L^2$-norm. Let $\tilde{T}_n$ be the associated "ideal" global test statistic, obtained by replacing $\hat{\boldsymbol{\beta}}(t)$ by the "ideal" estimator $\tilde{\boldsymbol{\beta}}(t)$ as defined in (3.3).

To derive the asymptotic random expression of $T_n$, we assume that $\gamma(s,t)$ has finite trace, that is, $\mathrm{tr}(\gamma) = \int_a^b \gamma(t,t)\,dt < \infty$. Let $\lambda_1, \lambda_2, \ldots$ be the eigenvalues, in decreasing order, and $\phi_1(t), \phi_2(t), \ldots$ be the associated orthonormal eigenfunctions of $\gamma(s,t)$. Let $m$ denote the number of positive eigenvalues. When all the eigenvalues are positive, we let $m = \infty$. Then $\lambda_r > 0$ for $r \leq m$ and $\lambda_r = 0$ for all $r > m$. Since $\mathrm{tr}(\gamma) < \infty$ implies $\int_a^b \int_a^b \gamma^2(s,t)\,ds\,dt < \infty$ by the Cauchy–Schwarz inequality, the covariance function $\gamma(s,t)$ has the singular value decomposition ([27], page 3)

$$(3.12) \qquad \gamma(s,t) = \sum_{r=1}^{m} \lambda_r \phi_r(s) \phi_r(t), \qquad s,t \in \mathcal{T} = [a,b].$$

THEOREM 7. *Assume the conditions of Theorem 6 are satisfied. Then as $n \to \infty$, we have*

$$(3.13) \qquad T_n = \tilde{T}_n + n^{1/2-(p+1)\delta/(2p+3)} O_P(1).$$

*In addition, assume $\delta > 1 + 1/[2(p+1)]$ and $\gamma(s,t)$ has finite trace so that it has the singular value decomposition (3.12). Then as $n \to \infty$, we have*

$$(3.14) \qquad T_n \stackrel{d}{=} \sum_{r=1}^{m} \lambda_r A_r + o_P(1), \qquad A_r \sim \chi_k^2(u_r^2),$$

*where $X \stackrel{d}{=} Y$ means the random variables $X$ and $Y$ have the same distribution, $\chi_k^2$ denotes a $\chi^2$-distribution with $k$ degrees of freedom and the noncentral parameters*

$$(3.15) \qquad u_r^2 = \lambda_r^{-1} \left\| \int_a^b \boldsymbol{\eta}_w(t) \phi_r(t)\,dt \right\|^2.$$

*Under $H_0$, $\boldsymbol{\eta}_w(t) \equiv 0$ so that all the $u_r^2$ are 0.*

Theorem 7 suggests that the distribution of $T_n$ is asymptotically the same as that of a $\chi^2$-type mixture. There are three possible methods that can be used to approximate the null distribution of $T_n$: $\chi^2$-approximation, simulation and bootstrapping. In the first two methods, we approximate the null distribution of $T_n$ by that of the $\chi^2$-type mixture $S = \sum_{r=1}^{\hat{m}} \hat{\lambda}_r A_r$, where $A_r \sim \chi_k^2$, $\hat{\lambda}_r$ are the eigenvalues of $\hat{\gamma}(s,t)$ and $\hat{m}$ is some well-chosen integer such that the eigenvalues $\hat{\lambda}_r, r = 1, 2, \ldots, \hat{m}$, explain a sufficiently large portion of the total variation $\mathrm{tr}(\hat{\gamma}) = \sum_{r=1}^{\infty} \hat{\lambda}_r$ and $\hat{\lambda}_r, r = \hat{m}+1, \hat{m}+2, \ldots$, are essentially 0. Besse [1] proposed a simple method for selecting such an $\hat{m}$. A simple and natural choice of $\hat{m}$ is the number of positive eigenvalues of $\hat{\gamma}(s,t)$. We found that the second method worked well in our simulation study and in the real data application presented in the next two sections.

In the $\chi^2$-approximation method, the distribution of $S$ is approximated by that of a random variable $R = \alpha\chi_d^2 + \beta$ via matching the first three cumulants of $R$ and $S$ to determine the unknown parameters $\alpha, d$ and $\beta$ [5, 29]. In the simulation method, the sampling distribution of $S$ is computed based on a sample of $S$ obtained via repeatedly generating $(A_1, A_2, \ldots, A_{\hat{m}})$. The bootstrap method is slightly more complicated. In the bootstrap method, we generate a sample of subject effects $v_i^*(t), i = 1, 2, \ldots, n$, from the estimated subject effects $\hat{v}_{i,1}(t), i = 1, 2, \ldots, n$, under $H_1$ and then construct a bootstrap sample, $f_i^*(t) = \mathbf{x}_i^T \hat{\boldsymbol{\beta}}_0(t) + v_i^*(t), i = 1, 2, \ldots, n$, where $\hat{\boldsymbol{\beta}}_0(t)$ is the estimator of $\boldsymbol{\beta}(t)$ under $H_0$ so that $\mathbf{C}\hat{\boldsymbol{\beta}}_0(t) = \mathbf{c}(t)$. Let $\hat{\boldsymbol{\beta}}^*(t)$ be the bootstrap estimator of $\boldsymbol{\beta}(t)$ based on the above bootstrap sample. We then use it to compute

$$T_n^* = \int_a^b \|\mathbf{w}^*(t)\|^2 \, dt = \sum_{l=1}^k \int_a^b w_l^{*2}(t) \, dt,$$

where $\mathbf{w}^*(t)$ can be obtained by replacing $\hat{\boldsymbol{\beta}}(t)$ with $\hat{\boldsymbol{\beta}}^*(t)$ in the definition (3.8) of $\mathbf{w}(t)$. The bootstrap null distribution of $T_n$ is obtained by the sampling distribution of $T_n^*$ via $B$ replications of the above bootstrap process for some large $B$, for example, $B = 10{,}000$.

**4. A simulation study.** In this section we aim to investigate the effect of the bandwidth selected by the GCV rule (2.16) on the average MSE (2.5) of the $p$-order LPK reconstructions $\hat{f}_i(t), i = 1, 2, \ldots, n$, and the MSE of the mean function estimator $\hat{\eta}(t)$ via a simulation study. We generated simulation samples from the model

$$\begin{aligned}
y_i(t) &= \eta(t) + v_i(t) + \varepsilon_i(t), \\
\eta(t) &= a_0 + a_1\phi_1(t) + a_2\phi_2(t), \\
v_i(t) &= b_{i0} + b_{i1}\psi_1(t) + b_{i2}\psi_2(t), \\
\mathbf{b}_i &= [b_{i0}, b_{i1}, b_{i2}]^T \sim N[0, \operatorname{diag}(\sigma_0^2, \sigma_1^2, \sigma_2^2)], \\
\varepsilon_i(t) &\sim N[0, \sigma_\varepsilon^2(1+t)], \qquad i = 1, 2, \ldots, n,
\end{aligned}$$

where $n$ is the number of subjects and $\mathbf{b}_i$ and $\varepsilon_i(t)$ are independent. The scheduled design time points are $t_j = j/(m+1), j = 1, 2, \ldots, m$. To obtain an unbalanced design which is more realistic, we randomly removed some responses on a subject at a rate $r_{miss}$ so that on average there are about $m(1-r_{miss})$ measurements on a subject, and $nm(1-r_{miss})$ measurements in a whole simulated sample. For simplicity, in this simulation the parameters we actually used are $[a_0, a_1, a_2] = [1.2, 2.3, 4.2], [\sigma_0^2, \sigma_1^2, \sigma_2^2, \sigma_\varepsilon^2] = [1, 2, 3, 0.1]$, $\phi_1(t) = \psi_1(t) = \cos(2\pi t), \phi_2(t) = \psi_2(t) = \sin(2\pi t), r_{miss} = 10\%$, $m = 40$ and $n = 20, 30$ and $40$.

For a simulated sample, the $p$-order LPK reconstructions $\hat{f}_i(t)$ were obtained using a local linear (i.e., $p=1$) smoother [8, 9] with the well-known Gaussian kernel. We considered five bandwidth choices, $0.5h^*, 0.8h^*, h^*, 1.25h^*$ and $2h^*$, where $h^*$ is the bandwidth selected by the GCV rule (2.16). For a simulated sample, the average MSE for $\hat{f}_i(t)$ and the MSE for the mean function estimator $\hat{\eta}(t)$ were computed respectively as $\text{MSE}_f = (nM)^{-1} \sum_{i=1}^{n} \sum_{j=1}^{M} \{\hat{f}_i(\tau_j) - f_i(\tau_j)\}^2$ and $\text{MSE}_\eta = M^{-1} \sum_{j=1}^{M} \{\hat{\eta}(\tau_j) - \eta(\tau_j)\}^2$, where $\tau_1, \ldots, \tau_M$ are $M$ time points equally-spaced in $[0, 1]$, for some large $M$, for example, $M = 400$.

Figure 1 presents the simulation results. The boxplots were based on 200 simulated samples. From left to right, panels are respectively for GCV, $\text{MSE}_f$ and $\text{MSE}_\eta$; from top to bottom, panels are respectively for $n = 20, 30$ and $40$. In each of the panels, the first five boxplots are associated with the five bandwidth choices: $0.5h^*, 0.8h^*, h^*, 1.25h^*$ and $2h^*$, respectively; the sixth boxplot in each of the $\text{MSE}_\eta$ panels is associated with the "ideal" estimator $\tilde{\eta}(t)$; see (2.10) for its definition.

From Figure 1, we may conclude that (a) overall, the GCV rule (2.16) performed well in the sense of choosing proper bandwidths to minimize the average MSE (2.5); (b) bandwidths smaller than $h^*$ help reduce the $\text{MSE}_\eta$ but do not by much, while bandwidths larger than $h^*$ do enlarge $\text{MSE}_\eta$ substantially; and (c) the $\text{MSE}_\eta$ based on $\hat{\eta}(t)$ and those based on the "ideal" estimator $\tilde{\eta}(t)$ are nearly the same unless the bandwidths are substantially larger than $h^*$.

**5. Application to the Canadian temperature data.** The Canadian temperature data (Canadian Climate Program [6]) were downloaded from <ftp://ego.psych.mcgill.ca/pub/ramsay/FDAfuns/Matlab/> at the book website of Ramsay and Silverman [23, 24]. The data are the daily temperature records of 35 Canadian weather stations over a year (365 days), among which 15 are in Eastern, another 15 in Western and the remaining five in Northern Canada. This is a typical functional data set with the number of measurements per subject ($n_i = 365$) being much larger than the number of subjects ($n = 35$). We shall use this functional data set only to illustrate the methodologies developed in this paper. For a more formal analysis, this functional data set should be first registered using either a parametric curve registration method proposed by Silverman [26] or a more flexible nonparametric curve registration method developed by Ramsay and Li [22]. Our methodologies can then be applied similarly to the resulting registered functional data set.

Figure 2 presents the individual curve reconstructions of the Canadian temperature data. These reconstructions were obtained by applying the local linear ($p = 1$) kernel fit [8, 9] with the well-known Gaussian kernel to the individual temperature records of each of the 35 weather stations, but with a common bandwidth $h^* = 2.79$, selected by the GCV rule (2.16). It can be seen that the Eastern weather station temperature curves (solid) mix up with the Western weather station temperature curves (dot-dashed), but most of the Eastern and Western weather station temperature curves stay higher than the Northern weather station temperature curves (dashed). This is reasonable since the Eastern and Western weather stations are located at about the same latitudes, while the Northern weather stations are located at higher latitudes.

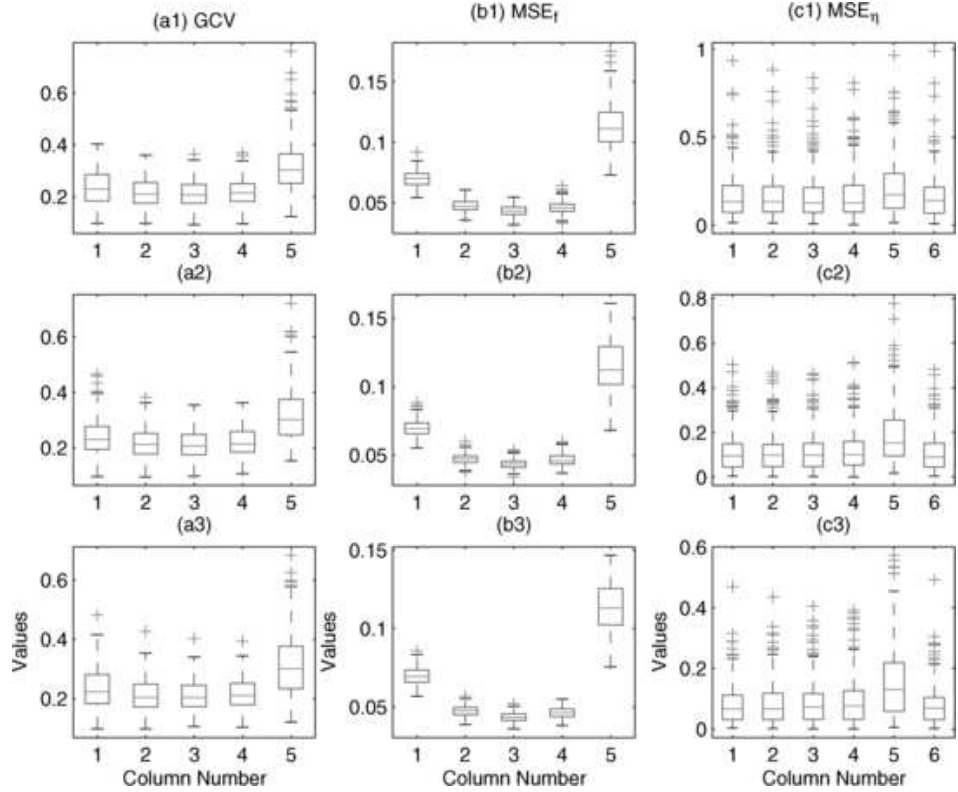

FIG. 1. *Simulation results. From left to right, panels are, respectively, for GCV, $MSE_f$ and $MSE_\eta$; from top to bottom, panels are, respectively, for $n = 20, 30$ and $40$. In each of the panels, the first five boxplots are associated with the five bandwidth choices $0.5h^*, 0.8h^*, h^*, 1.25h^*$ and $2h^*$, where $h^*$ is the GCV bandwidth; the sixth boxplot in a $MSE_\eta$ panel is associated with the "ideal" estimator $\tilde{\eta}(t)$.*

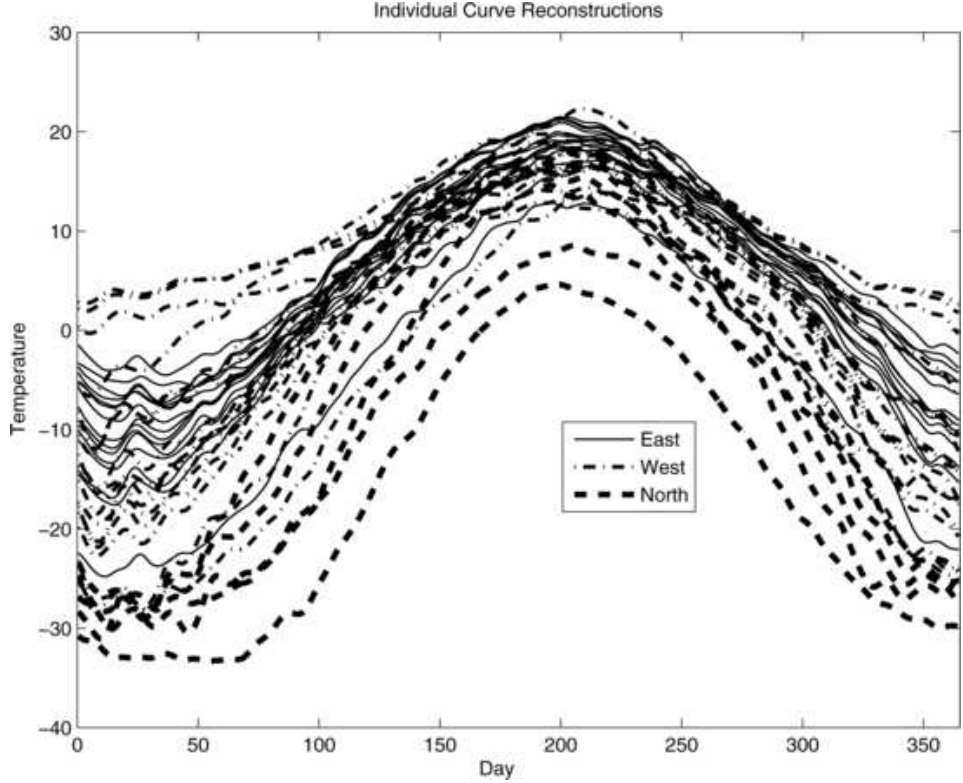

Fig. 2. *Local linear (p = 1) individual curve reconstructions of the Canadian temperature data with the bandwidth $h^* = 2.79$, selected by GCV. Eastern weather stations: solid curves; Western weather stations: dot-dashed curves; and Northern weather stations: dashed curves.*

We then modeled the Canadian temperature data set by the functional linear model (1.3) with the covariates

$$\mathbf{x}_i = \begin{cases} [1,0,0]^T, & \text{if weather station } i \text{ is located in Eastern Canada}, \\ [0,1,0]^T, & \text{if weather station } i \text{ is located in Western Canada}, \\ [0,0,1]^T, & \text{if weather station } i \text{ is located in Northern Canada}, \\ & i = 1, 2, \ldots, 35, \end{cases}$$

and the coefficient function vector $\boldsymbol{\beta}(t) = [\beta_1(t), \beta_2(t), \beta_3(t)]^T$, where $\beta_1(t)$, $\beta_2(t)$ and $\beta_3(t)$ are the covariate effect (mean temperature) functions of the Eastern, Western and Northern weather stations, respectively.

Figure 3 superimposes the estimated mean temperature functions of the Eastern, Western and Northern weather stations, together with their 95% standard deviation bands. Based on the 95% standard deviation bands, some informal conclusions can be made. First of all, over the whole year ($[a,b] = [1, 365]$), the differences between the mean temperature functions of the Eastern and the Western weather stations are much less significant than the differences between the mean temperature functions of the Eastern and the Northern weather stations, or between the Western and the Northern weather stations. This is because the 95% standard deviation band of the Eastern weather station mean temperature function covers (before Day 151) or stays close (after Day 151) to the mean temperature function of the Western weather stations; however, the 95% standard deviation bands of the Eastern and Western weather station mean temperature functions are far away from the mean temperature function of the Northern weather stations. Second, the significances of the differences between the mean temperature functions of the Eastern and the Western weather stations for different seasons are different. During the Spring (usually defined as the months of March, April and May or $[a,b] = [60, 151]$), the mean temperature functions are nearly the same, but this is not the case during the Summer (June, July and August or $[a,b] = [152, 243]$) or during the Autumn (September, October and November or $[a,b] = [244, 334]$). These conclusions can be made more clear via the hypothesis testing problem (3.6) with $t \in \mathcal{T} = [a, b]$ using the global testing statistic $T_n$ (3.11) and with $a, b, \mathbf{c}$ and $\mathbf{C}$ properly specified. For example, to test if the mean temperature functions of the Eastern and Western weather stations during the Spring are the same, we take $a = 60$, $b = 151$, $\mathbf{c} = 0$ and $\mathbf{C} = [1, -1, 0]$; and to test if the mean temperature functions of the Eastern, Western and Northern weather stations during the Autumn are the same, we take $a = 244$, $b = 334$, $\mathbf{c} = [0, 0]^T$ and

$$\mathbf{C} = \begin{bmatrix} 1, & 0, & -1 \\ 0, & 1, & -1 \end{bmatrix}.$$

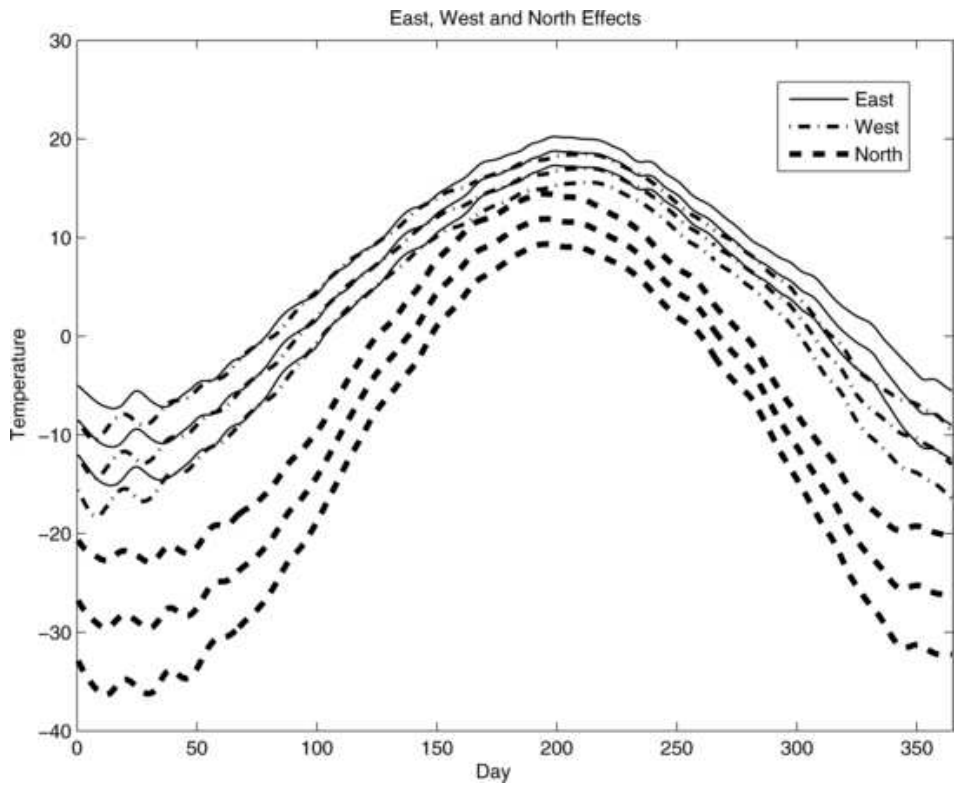

Fig. 3. *Estimated mean temperature functions of the Eastern, Western and Northern weather stations with 95% standard deviation bands (Eastern weather stations, solid; Western weather stations, dot-dashed; and Northern weather stations, dashed).*

TABLE 1
*Significance test results for the differences of the mean temperature functions of the Eastern and Western weather stations based on* 10,000 *replications*

| [a, b] | $h$ | $T_n$ | P-values $\chi^2$-approximation | Simulation | Bootstrapping |
|---|---|---|---|---|---|
| [1, 365] | $h^*/2$ | 59954 | 0.179 | 0.179 | 0.166 |
| (Whole year) | $h^*$ | 58248 | 0.185 | 0.181 | 0.180 |
| | $2h^*$ | 56868 | 0.189 | 0.185 | 0.184 |
| [60, 151] | $h^*/2$ | 945 | 0.842 | 0.836 | 0.834 |
| (Spring) | $h^*$ | 656 | 0.940 | 0.874 | 0.877 |
| | $2h^*$ | 378 | 1.000 | 0.923 | 0.922 |
| [152, 243] | $h^*/2$ | 6625 | 0.078 | 0.075 | 0.068 |
| (Summer) | $h^*$ | 6432 | 0.082 | 0.084 | 0.083 |
| | $2h^*$ | 6322 | 0.085 | 0.086 | 0.075 |
| [244, 334] | $h^*/2$ | 28748 | 0.011 | 0.011 | 0.009 |
| (Autumn) | $h^*$ | 28303 | 0.012 | 0.013 | 0.008 |
| | $2h^*$ | 27526 | 0.014 | 0.015 | 0.010 |

We first tested the differences of the mean temperature functions of the Eastern and Western Canadian weather stations for the whole year, and during the Spring, Summer and Autumn. Table 1 shows the significance test results, where the simulation and bootstrap P-values were computed based on 10,000 replications. For each choice of the seasonal period $[a, b]$, we used three different bandwidth choices, $h^*/2$, $h^*$ and $2h^*$, where $h^* = 2.79$ was selected by the GCV rule (2.16). For each bandwidth choice, the associated test statistics $T_n$ were computed using (3.11). For each $T_n$, we computed its P-value using the $\chi^2$-approximation, simulation and bootstrap methods which were described briefly in Section 3.2. Figure 4 displays the null probability density function (p.d.f.) approximations obtained using the three methods. It seems that all three approximations perform reasonably well except at the left boundary where the $\chi^2$-approximations seem problematic. Nevertheless, from the table, we can see that the significance test results are not strongly affected by the bandwidths used; moreover, we can see that the differences between the mean temperature functions of the Eastern and Western weather stations over the whole year (P-value $\geq 0.166$) are larger than their differences during the Spring (P-value $\geq 0.834$), but much smaller than their differences during the Summer (P-value $< 0.068$) or during the Autumn (P-value $< 0.015$). These results are consistent with those observed from Figure 3.

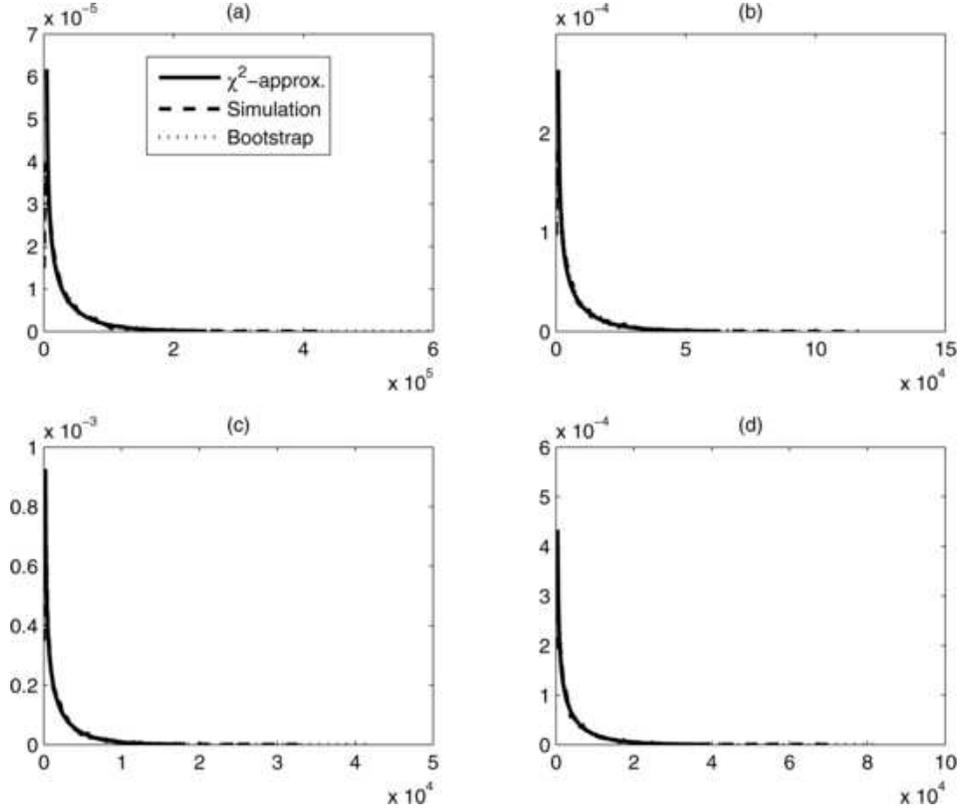

Fig. 4. *Null p.d.f. approximations ($\chi^2$-approximation, solid; simulation, dashed; bootstrap, dotted) of the global test statistic $T_n$ (3.11) when $h^* = 2.79$. (a) $[a,b] = [1, 365]$; (b) $[a,b] = [60, 151]$; (c) $[a,b] = [152, 243]$; and (d) $[a,b] = [244, 334]$.*

Following the same procedure, we also tested the following null hypotheses: the mean temperature functions are the same between (1) the Eastern and Northern; (2) the Western and Northern; and (3) the Eastern, Western and Northern weather stations for the following periods: (1) the whole year; (2) the Spring; (3) the Summer; and (4) the Autumn. As expected, we rejected all these null hypotheses with P-value 0. These results are also consistent with those observed from Figure 3.

**6. Technical proofs.** In this section we outline the technical proofs of some of the asymptotic results. Before we proceed, we list the following useful lemmas. Proof of the first lemma can be found in [10], page 64. Notice that, under Condition A4, "$n \to \infty$" implies that "$n_i \to \infty$."

LEMMA 1. *Assume Condition A is satisfied. Then as $n \to \infty$, we have*

$$K_h^{n_i}(t_{ij} - t) = \frac{1}{n_i \pi(t)} K_h^*(t_{ij} - t)[1 + o_P(1)],$$

*where $K^*(\cdot)$ is the LPK equivalent kernel ([10], page 64).*

LEMMA 2. *We always have*

$$\sum_{j=1}^{n_i} K_h^{n_i}(t_{ij} - t)(t_{ij} - t)^r = \begin{cases} 1, & \text{when } r = 0, \\ 0, & \text{otherwise.} \end{cases}$$

*Assume Condition A is satisfied. Then as $n \to \infty$, we have*

$$\sum_{j=1}^{n_i} K_h^{n_i}(t_{ij} - t)(t_{ij} - t)^{p+1} = \frac{B_{p+1}(K^*)h^{p+1}}{\pi(t)}[1 + o_P(1)],$$

$$\sum_{j=1}^{n_i} \{K_h^{n_i}(t_{ij} - t)\}^2 = \frac{V(K^*)}{\pi(t)}(n_i h)^{-1}[1 + o_P(1)],$$

$$\sum_{j=1}^{n_i} K_h^{n_i}(t_{ij} - s)K_h^{n_i}(t_{ij} - t) = \frac{K^{*(1)}((s-t)/h)}{\pi(t)}(n_i h)^{-1}[1 + o_P(1)],$$

*where $B_r(\cdot)$ and $V(\cdot)$ are defined in (2.4).*

Let $r_i(t) = \hat{f}_i(t) - f_i(t), i = 1, 2, \ldots, n$, where $\hat{f}_i(t)$ are the $p$-order LPK reconstructions of $f_i(t)$ given in Section 2.1. Let $\bar{r}(t) = n^{-1} \sum_{i=1}^n r_i(t)$ and $\bar{f}(t) = n^{-1} \sum_{i=1}^n f_i(t)$. Using Lemmas 1 and 2, we can prove the following useful lemma.

LEMMA 3. *Assume Condition A is satisfied. Then as $n \to \infty$, we have*

$$E\{r_i(t)|\mathcal{D}\} = \frac{B_{p+1}(K^*)\eta^{(p+1)}(t)}{(p+1)!} h^{p+1}[1 + o_P(1)],$$

$$\mathrm{Cov}\{r_i(s), r_i(t)|\mathcal{D}\} = \left\{ \frac{K^{*(1)}((s-t)/h)}{\pi(t)}(n_i h)^{-1} \right.$$

$$\left. + \frac{B_{p+1}^2(K^*)\gamma_{p+1,p+1}(s,t)}{(p+1)!^2} h^{2(p+1)} \right\}[1 + o_P(1)],$$

$$\mathrm{Cov}\{r_i(s), f_i(t)|\mathcal{D}\} = \frac{B_{p+1}(K^*)\gamma_{p+1,0}(s,t)}{(p+1)!} h^{p+1}[1 + o_P(1)].$$

PROOF. By (2.3) and Lemma 1, we have

$$r_i(t) = \sum_{j=1}^{n_i} K_h^{n_i}(t_{ij} - t)\varepsilon_{ij} + \sum_{j=1}^{n_i} K_h^{n_i}(t_{ij} - t)\{f_i(t_{ij}) - f_i(t)\}.$$

It follows that

$$\mathrm{E}(r_i(t)|\mathcal{D}) = \sum_{j=1}^{n_i} K_h^{n_i}(t_{ij} - t)\{\eta(t_{ij}) - \eta(t)\}.$$

Applying Taylor's expansion and Lemmas 1 and 2, we have

(6.1)
$$\mathrm{E}(r_i(t)|\mathcal{D}) = \sum_{j=1}^{n_i} K_h^{n_i}(t_{ij} - t)\left\{\sum_{l=1}^{p+1} \eta^{(l)}(t)\frac{(t_{ij}-t)^l}{l!} + o[(t_{ij}-t)^{p+1}]\right\}$$

$$= \frac{B_{p+1}(K^*)\eta^{(p+1)}(t)}{(p+1)!} h^{p+1}[1 + o_P(1)].$$

Similarly, by the independence of $f_i(t)$ and $\varepsilon_i(t)$, we have

$$\mathrm{Cov}(r_i(s), r_i(t)|\mathcal{D})$$

$$= \sum_{j=1}^{n_i} K_h^{n_i}(t_{ij} - s)K_h^{n_i}(t_{ij} - t)\sigma^2(t_{ij})$$

$$+ \sum_{j=1}^{n_i}\sum_{l=1}^{n_i} K_h^{n_i}(t_{ij} - s)K_h^{n_i}(t_{il} - t)$$

(6.2)
$$\times \{\gamma(t_{ij}, t_{il}) - \gamma(t_{ij}, t) - \gamma(s, t_{il}) + \gamma(s, t)\}$$

$$= \left\{ \frac{K^{*(1)}[(s-t)/h]\sigma^2(s)}{\pi(t)}(n_i h)^{-1} \right.$$

$$\left. + \frac{B_{p+1}^2(K^*)\gamma_{p+1,p+1}(s,t)}{(p+1)!^2} h^{2(p+1)} \right\}[1 + o_P(1)].$$

In particular, letting $s = t$, we obtain

$$\mathrm{Var}(r_i(t)|\mathcal{D}) = \left\{ \frac{V(K^*)\sigma^2(t)}{\pi(t)}(n_i h)^{-1} \right.$$

(6.3)
$$+\frac{B_{p+1}^2(K^*)\gamma_{p+1,p+1}(t,t)}{(p+1)!^2}h^{2(p+1)}\bigg\}[1+o_P(1)],$$

as desired. Lemma 3 is proved. $\square$

Direct application of Lemma 3 leads to the following.

LEMMA 4. *Assume Condition A is satisfied. Then as $n \to \infty$, we have*

$$\mathrm{E}(\bar{r}(t)|\mathcal{D}) = \frac{B_{p+1}(K^*)\eta^{(p+1)}(t)}{(p+1)!}h^{p+1}[1+o_P(1)],$$

$$\mathrm{Cov}(\bar{r}(s),\bar{r}(t)|\mathcal{D}) = n^{-1}\bigg\{\frac{K^{*(1)}((s-t)/h)}{\pi(t)}(\tilde{m}h)^{-1}$$
$$+\frac{B_{p+1}^2(K^*)\gamma_{p+1,p+1}(s,t)}{(p+1)!^2}h^{2(p+1)}\bigg\}[1+o_P(1)],$$

$$\mathrm{Cov}(\bar{r}(s),\bar{f}(t)|\mathcal{D}) = \frac{B_{p+1}(K^*)\gamma_{p+1,0}(s,t)}{(p+1)!}n^{-1}h^{p+1}[1+o_P(1)],$$

*where $\tilde{m} = (n^{-1}\sum_{i=1}^n n_i^{-1})^{-1}$, as defined in Theorem 2.*

PROOF OF THEOREM 1. For each $i=1,2,\ldots,n$, by (6.1) and (6.3), we have

(6.4)
$$\mathrm{E}\{[\hat{f}_i(t)-f_i(t)]^2|\mathcal{D}\}$$
$$= \mathrm{E}\{r_i^2(t)|\mathcal{D}\} = \{\mathrm{E}(r_i(t)|\mathcal{D})\}^2 + \mathrm{Var}(r_i(t)|\mathcal{D})$$
$$= \bigg\{\frac{B_{p+1}^2(K^*)[(\eta^{(p+1)}(t))^2+\gamma_{p+1,p+1}(t,t)]}{(p+1)!^2}h^{2(p+1)}$$
$$+\frac{V(K^*)\sigma^2(t)}{\pi(t)}(n_ih)^{-1}\bigg\}[1+o_P(1)].$$

Theorem 1, that is, the expression (2.5), then follows directly. $\square$

PROOF OF THEOREM 2. Under Condition A, the coefficients of $h^{2(p+1)}$ and $(n_ih)^{-1}$ in the expression (6.4) are uniformly bounded over the finite interval $\mathcal{T}=[a,b]$. Moreover, since $n_i \geq Cn^\delta$ and $h = O(n^{-\delta/(2p+3)})$, we have $O(h^{2(p+1)}) = O((n_ih)^{-1}) = O(n^{-2(p+1)\delta/(2p+3)}) = n^{-2(p+1)\delta/(2p+3)}O(1)$. Thus, $\mathrm{E}\{r_i^2(t)|\mathcal{D}\} = O_{\mathrm{UP}}[h^{2(p+1)}+(n_ih)^{-1}] = n^{-2(p+1)\delta/(2p+3)}O_{\mathrm{UP}}(1)$. Therefore, $\hat{f}_i(t) = f_i(t) + n^{-(p+1)\delta/(2p+3)}O_{\mathrm{UP}}(1)$. Theorem 2 is then proved. $\square$

PROOF OF THEOREM 3. First of all, notice that $\hat{\eta}(t) = n^{-1}\sum_{i=1}^n \hat{f}_i(t) = \bar{f}(t) + \bar{r}(t)$. It follows that $\mathrm{Bias}(\hat{\eta}(t)|\mathcal{D}) = \mathrm{E}(\bar{r}(t)|\mathcal{D})$, $\mathrm{Cov}(\hat{\eta}(s),\hat{\eta}(t)|\mathcal{D}) =$ $\operatorname{Cov}(\bar{f}(s), \bar{f}(t)) + \operatorname{Cov}(\bar{f}(s), \bar{r}(t)) + \operatorname{Cov}(\bar{r}(s), \bar{f}(t)) + \operatorname{Cov}(\bar{r}(s), \bar{r}(t))$. The results of Theorem 3 follow directly from Lemma 4. □

PROOF OF THEOREM 4. Since $\hat{\eta}(t) = \bar{f}(t) + \bar{r}(t) = \tilde{\eta}(t) + \bar{r}(t)$, in order to show the first expression in (2.11), it is sufficient to prove that $\operatorname{E}\{\bar{r}^2(t)|\mathcal{D}\} = n^{-2(p+1)\delta/(2p+3)} O_{\mathrm{UP}}(1)$. This result follows directly from $\operatorname{E}\{\bar{r}^2(t)|\mathcal{D}\} = \{\operatorname{E}(\bar{r}(t)|\mathcal{D})\}^2 + \operatorname{Var}(\bar{r}(t)|\mathcal{D})$ and Lemma 4. To show the second expression in (2.11), notice that the covariance estimator $\hat{\gamma}(s,t)$ can be expressed as

$$\hat{\gamma}(s,t) = \frac{1}{n} \sum_{i=1}^{n} \{f_i(s) - \bar{f}(s)\}\{f_i(t) - \bar{f}(t)\}$$

$$+ \frac{1}{n} \sum_{i=1}^{n} \{f_i(s) - \bar{f}(s)\}\{r_i(t) - \bar{r}(t)\}$$

$$+ \frac{1}{n} \sum_{i=1}^{n} \{r_i(s) - \bar{r}(s)\}\{f_i(t) - \bar{f}(t)\}$$

$$+ \frac{1}{n} \sum_{i=1}^{n} \{r_i(s) - \bar{r}(s)\}\{r_i(t) - \bar{r}(t)\}$$

$$\equiv \tilde{\gamma}(s,t) + I_1 + I_2 + I_3,$$

where $r_i(t) = \hat{f}_i(t) - f_i(t), i = 1, 2, \ldots, n$, are independent and asymptotically have the same variance. By the law of large numbers and by Lemma 3, we have

$$I_1 = \operatorname{E}\left\{ n^{-1} \sum_{i=1}^{n} \operatorname{E}[(f_i(s) - \bar{f}(s))(r_i(t) - \bar{r}(t))|\mathcal{D}] \right\} O_P(1)$$

$$= \operatorname{E}\{\operatorname{Cov}(f_1(s), r_1(t)|\mathcal{D})\} O_P(1)$$

$$= n^{-(p+1)\delta/(2p+3)} O_{\mathrm{UP}}(1).$$

Similarly, we can show that

$$I_2 = n^{-(p+1)\delta/(2p+3)} O_{\mathrm{UP}}(1) \quad \text{and} \quad I_3 = n^{-2(p+1)\delta/(2p+3)} O_{\mathrm{UP}}(1).$$

The second expression in (2.11) then follows. When $\delta > 1 + 1/[2(p+1)]$, we have

$$n^{1/2}\{\tilde{\eta}(t) - \hat{\eta}(t)\} = o_{\mathrm{UP}}(1), \qquad n^{1/2}\{\tilde{\gamma}(s,t) - \hat{\gamma}(s,t)\} = o_{\mathrm{UP}}(1).$$

By the definition of $\tilde{\eta}(t)$ and $\tilde{\gamma}(s,t)$, we have

$$\tilde{\eta}(t) = \eta(t) + \bar{v}(t), \qquad \tilde{\gamma}(s,t) = n^{-1} \sum_{i=1}^{n} v_i(s) v_i(t) - \bar{v}(s) \bar{v}(t).$$

By the law of large numbers and the central limit theorem, it is easy to show that

$$n^{1/2}\{\tilde{\eta}(t) - \eta(t)\} \sim \text{AGP}(0,\gamma), \qquad n^{1/2}\{\tilde{\gamma}(s,t) - \gamma(s,t)\} \sim \text{AGP}(0,\gamma^*),$$

where

$$\gamma^*\{(s_1,t_1),(s_2,t_2)\} = \text{Cov}\{v_1(s_1)v_1(t_1), v_1(s_2)v_1(t_2)\}$$
$$= \text{E}\{v_1(s_1)v_1(t_1)v_1(s_2)v_1(t_2)\} - \gamma(s_1,t_1)\gamma(s_2,t_2).$$

In particular, when $v(t)$ is a Gaussian process, we have

$$\text{E}\{v_1(s_1)v_1(t_1)v_1(s_2)v_1(t_2)\}$$
$$= \gamma(s_1,t_1)\gamma(s_2,t_2) + \gamma(s_1,t_2)\gamma(s_2,t_1) + \gamma(s_1,s_2)\gamma(t_1,t_2).$$

Thus, $\gamma^*\{(s_1,t_1),(s_2,t_2)\} = \gamma(s_1,t_2)\gamma(s_2,t_1) + \gamma(s_1,s_2)\gamma(t_1,t_2)$. The proof of Theorem 4 is finished. □

PROOF OF THEOREM 5. Under Condition A and by Theorem 2, we have $\hat{f}_i(t_{ij}) = f_i(t_{ij}) + n^{-(p+1)\delta/(2p+3)}O_{\text{UP}}(1)$. It follows that

$$\hat{\varepsilon}_{ij}^2 = \{y_{ij} - \hat{f}_i(t_{ij})\}^2 = \{\varepsilon_{ij} + n^{-(p+1)\delta/(2p+3)}O_{\text{UP}}(1))\}^2$$
$$= \varepsilon_{ij}^2 + 2n^{-(p+1)\delta/(2p+3)}\varepsilon_{ij}O_{\text{UP}}(1) + n^{-2(p+1)\delta/(2p+3)}O_{\text{UP}}(1).$$

Plugging this into (2.14) with $b = O(N^{-1/5})$, we have $\hat{\sigma}^2(t) = I_1 + I_2 + I_3$, where under the given conditions and by standard kernel estimation theory,

$$I_1 = \frac{\sum_{i=1}^n \sum_{j=1}^{n_i} H_b(t_{ij} - t)\varepsilon_{ij}^2}{\sum_{i=1}^n \sum_{j=1}^{n_i} H_b(t_{ij} - t)} = \sigma^2(t) + N^{-2/5}O_{\text{UP}}(1),$$

$$I_2 = 2\frac{\sum_{i=1}^n \sum_{j=1}^{n_i} H_b(t_{ij} - t)n^{-(p+1)\delta/(2p+3)}\varepsilon_{ij}O_{\text{UP}}(1)}{\sum_{i=1}^n \sum_{j=1}^{n_i} H_b(t_{ij} - t)}$$
$$= n^{-(p+1)\delta/(2p+3)}O_{\text{UP}}(1),$$

$$I_3 = \frac{\sum_{i=1}^n \sum_{j=1}^{n_i} H_b(t_{ij} - t)n^{-2(p+1)\delta/(2p+3)}O_{\text{UP}}(1)}{\sum_{i=1}^n \sum_{j=1}^{n_i} H_b(t_{ij} - t)}$$
$$= n^{-2(p+1)\delta/(2p+3)}O_{\text{UP}}(1).$$

Under Condition A4, $n_i \geq Cn^\delta$. This implies that $N = \sum_{i=1}^n n_i > Cn^{1+\delta}$. Thus, $N^{-2/5} = O(n^{-2(1+\delta)/5})$. It follows that $\hat{\sigma}^2(t) = \sigma^2(t) + O_{\text{UP}}(n^{-2(1+\delta)/5} + n^{-(p+1)\delta/(2p+3)})$, as desired. The proof of the theorem is completed. □

PROOF OF THEOREM 6. Under the conditions of Theorem 2, we have $|r_i(t)| = |\hat{f}_i(t) - f_i(t)| \leq n^{-(p+1)\delta/(2p+3)}C$ for some $C > 0$ for all $i$ and $t$. Let

$\Delta(t) = [\Delta_1(t), \ldots, \Delta_q(t)]^T = (\mathbf{X}^T\mathbf{X})^{-1}\mathbf{X}^T(\hat{\mathbf{f}}(t) - \mathbf{f}(t))$. Then for $r = 1, 2, \ldots, q$, we have

$$|\Delta_r(t)| = \left| n^{-1} \sum_{i=1}^n \mathbf{e}_{r,q}^T \left( n^{-1} \sum_{j=1}^n \mathbf{x}_j \mathbf{x}_j^T \right)^{-1} \mathbf{x}_i r_i(t) \right|$$

$$\leq n^{-1} \sum_{i=1}^n \left| \mathbf{e}_{r,q}^T \left( n^{-1} \sum_{j=1}^n \mathbf{x}_j \mathbf{x}_j^T \right)^{-1} \mathbf{x}_i \right| |r_i(t)|$$

$$\leq C n^{-(p+1)\delta/(2p+3)} \mathrm{E} |\mathbf{e}_{r,q}^T \Omega^{-1} \mathbf{x}_1| [1 + o_p(1)].$$

It follows that $\Delta(t) = n^{-(p+1)\delta/(2p+3)} O_{\mathrm{UP}}(1)$. The first expression in (3.4) follows directly from the fact $\hat{\boldsymbol{\beta}}(t) - \tilde{\boldsymbol{\beta}}(t) = \Delta(t)$.

To show the second expression in (3.4), notice that $\hat{v}_i(t) = \tilde{v}_i(t) + r_i(t) + \mathbf{x}_i^T[\hat{\boldsymbol{\beta}}(t) - \tilde{\boldsymbol{\beta}}(t)] = \tilde{v}_i(t) + n^{-(p+1)\delta/(2p+3)} O_{\mathrm{UP}}(1)$ because under the given conditions, we have $\mathbf{x}_i = O_{\mathrm{UP}}(1)$, $r_i(t) = n^{-(p+1)\delta/(2p+3)} O_{\mathrm{UP}}(1)$, and $\hat{\boldsymbol{\beta}}(t) - \tilde{\boldsymbol{\beta}}(t) = n^{-(p+1)\delta/(2p+3)} O_{\mathrm{UP}}(1)$. Further, by Condition B, we have $v_i(t) = O_{\mathrm{UP}}(1)$, therefore, $\hat{v}_i(s)\hat{v}_i(t) = \tilde{v}_i(s)\tilde{v}_i(t) + n^{-(p+1)\delta/(2p+3)} O_{\mathrm{UP}}(1)$. The second expression in (3.4) follows immediately.

When $\delta > 1 + 1/[2(p+1)]$, we have $(p+1)\delta/(2p+3) > 1/2$. Therefore, $\sqrt{n}[\hat{\boldsymbol{\beta}}(t) - \tilde{\boldsymbol{\beta}}(t)] = n^{1/2 - (p+1)\delta/(2p+3)} O_{\mathrm{UP}}(1) = o_{\mathrm{UP}}(1)$. Moreover, it is easy to show that

(6.5) $$\sqrt{n}[\tilde{\boldsymbol{\beta}}(t) - \boldsymbol{\beta}(t)] \sim \mathrm{AGP}(0, \gamma_\beta),$$

where $\gamma_\beta(s, t) = \gamma(s, t)\boldsymbol{\Omega}^{-1}$. The result in (3.5) follows immediately. The proof of the theorem is completed. $\square$

PROOF OF THEOREM 7. Recall that $\mathbf{w}(t) = [\mathbf{C}(\mathbf{X}^T\mathbf{X})^{-1}\mathbf{C}^T]^{-1/2}[\mathbf{C}\hat{\boldsymbol{\beta}}(t) - \mathbf{c}(t)]$, as defined in (3.8). Define $\tilde{\mathbf{w}}(t)$ similarly by replacing $\hat{\boldsymbol{\beta}}(t)$ with $\tilde{\boldsymbol{\beta}}(t)$. Then by (3.11), we have $T_n = \int_a^b \|\mathbf{w}(t)\|^2 \, dt$ and similarly, $\tilde{T}_n = \int_a^b \|\tilde{\mathbf{w}}(t)\|^2 \, dt$.

Let $\Delta(t) = \mathbf{w}(t) - \tilde{\mathbf{w}}(t) = [\mathbf{C}(\mathbf{X}^T\mathbf{X})^{-1}\mathbf{C}^T]^{-1/2}\mathbf{C}[\hat{\boldsymbol{\beta}}(t) - \tilde{\boldsymbol{\beta}}(t)]$. Then under the given conditions and by Theorem 6, we can show that $\Delta(t) = n^{1/2 - (p+1)\delta/(2p+3)} \times O_{\mathrm{UP}}(1)$. It follows that $\mathbf{w}(t) = \tilde{\mathbf{w}}(t) + n^{1/2 - (p+1)\delta/(2p+3)} O_{\mathrm{UP}}(1)$ and, hence, $T_n = \tilde{T}_n + 2\int_a^b \tilde{\mathbf{w}}(t)^T \Delta(t) \, dt + \int_a^b \|\Delta(t)\|^2 \, dt = \tilde{T}_n + n^{1/2 - (p+1)\delta/(2p+3)} O_p(1)$, as desired.

When $\delta > 1 + 1/[2(p+1)]$, we have $T_n = \tilde{T}_n + o_P(1)$ as $n \to \infty$. Thus, to show (3.14), it is sufficient to show $\tilde{T}_n \stackrel{\mathrm{d}}{=} \sum_{r=1}^m \lambda_r A_r + o_P(1)$. Using (6.5) in the proof of Theorem 6 above, it is easy to show that $\tilde{\mathbf{w}}(t) \sim \mathrm{AGP}(\boldsymbol{\eta}_w, \boldsymbol{\gamma}_w)$, where $\boldsymbol{\eta}_w(t) = \sqrt{n}(\mathbf{C}\Omega^{-1}\mathbf{C}^T)^{-1/2}[\mathbf{C}\boldsymbol{\beta}(t) - \mathbf{c}(t)]$ and $\boldsymbol{\gamma}_w(s, t) = \gamma(s, t)\mathbf{I}_k$, as defined in (3.9). It follows that the $k$ components of $\tilde{\mathbf{w}}(t)$ are independent of each other, and the $l$th component $\tilde{w}_l(t) \sim \mathrm{AGP}(\eta_{wl}, \gamma)$, where $\eta_{wl}(t)$ is the $l$th component of $\boldsymbol{\eta}_w(t)$ as defined in (3.9). Since $\gamma(s,t)$ has the singular value decomposition (3.12), we have $\tilde{w}_l(t) = \sum_{r=1}^{m} \xi_{lr} \phi_r(t)$, where

$$(6.6) \qquad \xi_{lr} = \int_a^b \tilde{w}_l(t) \phi_r(t) \, dt \sim \mathrm{AN}(\mu_{lr}, \lambda_r),$$

with $\mu_{lr} = \int_a^b \eta_{wl}(t) \phi_r(t) \, dt$. It follows that

$$\tilde{T}_n = \int_a^b \|\tilde{\mathbf{w}}(t)\|^2 \, dt = \sum_{l=1}^{k} \int_a^b \tilde{w}_l^2(t) \, dt$$
$$= \sum_{l=1}^{k} \sum_{r=1}^{m} \xi_{lr}^2 = \sum_{r=1}^{m} \sum_{l=1}^{k} \xi_{lr}^2$$

because the eigenfunctions $\phi_r(t)$ are orthonormal over $\mathcal{T} = [a,b]$ and the summation is exchangeable due to the nonnegativity of $\xi_{lr}^2$. By (6.6), we have $\sum_{l=1}^{k} \xi_{lr}^2 \stackrel{\mathrm{d}}{=} \lambda_r A_r$, where $A_r \sim \chi_k^2(u_r^2)$ with $u_r^2 = \lambda_r^{-1} \sum_{l=1}^{k} \mu_{lr}^2 = \lambda_r^{-1} \| \int_a^b \boldsymbol{\eta}_w(t) \times \phi_r(t) \, dt \|^2$, as given in (3.15). It follows that $\tilde{T}_n \stackrel{\mathrm{d}}{=} \sum_{r=1}^{m} \lambda_r A_r + o_P(1)$, as desired. The proof of the theorem is completed. $\square$

**Acknowledgments.** The authors thank the Editor, the Associate Editor and two reviewers for their helpful comments and invaluable suggestions that helped improve the paper substantially.

DEPARTMENT OF STATISTICS
AND APPLIED PROBABILITY
NATIONAL UNIVERSITY OF SINGAPORE
LOWER KENT RIDGE ROAD
3 SCIENCE DRIVE 2, SINGAPORE 119260
SINGAPORE
E-MAIL: stazjt@nus.edu.sg

DEPARTMENT OF MATHEMATICS
AND STATISTICS
SAN DIEGO STATE UNIVERSITY
5500 CAMPANILE DRIVE
SAN DIEGO, CA 92182
USA
E-MAIL: jchen@sciences.sdsu.edu